\documentclass[11pt]{article}
\usepackage{graphicx}
\usepackage{color}
\usepackage{amsmath}
\usepackage{amssymb}
\usepackage{amscd}
\usepackage{framed}
\usepackage{bbm}

\usepackage{fancyhdr,a4wide}
\usepackage{amsthm}

\newcommand{\R}{\mathbb{R}}
\newcommand{\inr}[1]{\left\langle #1 \right\rangle}

\newcommand{\E}{\mathbb{E}}

\newcommand{\eps}{\varepsilon}

\newtheorem{Theorem}{Theorem}[section]
\newtheorem{Lemma}[Theorem]{Lemma}

\newtheorem{Definition}[Theorem]{Definition}

\newtheorem{Corollary}[Theorem]{Corollary}
\newtheorem{Remark}[Theorem]{Remark}
\newtheorem{Example}[Theorem]{Example}

\numberwithin{equation}{section}

\def \proof {\noindent {\bf Proof.}\ \ }

\def \endproof
{{\mbox{}\nolinebreak\hfill\rule{2mm}{2mm}\par\medbreak}}
\def\IND{\mathbbm{1}}

\def\IND{\mathbbm{1}}

\begin{document}
\title{Column randomization and almost-isometric embeddings}
\author{Shahar Mendelson
\thanks{Centre for Mathematics and its Applications, The Australian National University. Email: shahar.mendelson@anu.edu.au}
}

\maketitle

\begin{abstract}
The matrix $A:\R^n \to \R^m$ is $(\delta,k)$-regular if for any $k$-sparse vector $x$,
$$
\left| \|Ax\|_2^2-\|x\|_2^2\right| \leq \delta \sqrt{k} \|x\|_2^2.
$$
We show that if $A$ is $(\delta,k)$-regular for $1 \leq k \leq 1/\delta^2$, then by multiplying the columns of $A$ by independent random signs, the resulting random ensemble $A_\eps$ acts on an arbitrary subset $T \subset \R^n$ (almost) as if it were gaussian, and with the optimal probability estimate: if $\ell_*(T)$ is the gaussian mean-width of $T$ and $d_T=\sup_{t \in T} \|t\|_2$, then with probability at least $1-2\exp(-c(\ell_*(T)/d_T)^2)$,
$$
\sup_{t \in T} \left| \|A_\eps t\|_2^2-\|t\|_2^2 \right| \leq C\left(\Lambda d_T \delta\ell_*(T)+(\delta \ell_*(T))^2 \right),
$$
where $\Lambda=\max\{1,\delta^2\log(n\delta^2)\}$. This estimate is optimal for $0<\delta \leq 1/\sqrt{\log n}$.
\end{abstract}

\section{Introduction}
Linear operators that act in an almost-isometric way on subsets of $\R^n$ are of obvious importance. Although approximations of isometries are the only operators that almost preserve the Euclidean norm of any point in $\R^n$, one may consider a more flexible alternative: a random ensemble of operators $\Gamma$ such that, for any fixed $T \subset \R^n$, with high probability, $\Gamma$ ``acts well" on every element of $T$. Such random ensembles have been studied extensively over the years, following the path paved by the celebrated work of Johnson and Lindenstrauss in \cite{MR737400}. Here we formulate the Johnson-Lindenstrauss Lemma in one of its gaussian versions:
\begin{Theorem} \label{thm:JL}
There exist absolute constants $c_0$ and $c_1$ such that the following holds. Let $1 \leq m \leq n$ and set $\Gamma:\R^n \to \R^m$ to be a random matrix whose entries are independent, standard gaussian random variables. Let $T \subset S^{n-1}$ be of cardinality at most $\exp(c_0m)$. Then for $m^{-1/2}\sqrt{\log |T|}<\rho<1$, with probability at least $1-2\exp(-c_1 \rho^2 m)$, for every $t \in T$,
$$
\left| \left\|m^{-1/2} \Gamma t\right\|_2^2 - 1 \right| \leq \rho.
$$
\end{Theorem}

The scope of Theorem \ref{thm:JL} can be extended to more general random ensembles than the gaussian one, e.g., to a random matrix whose rows are iid copies of a centred random vector that exhibits suitable decay properties (see, e.g. \cite{MR3552848,MR3565471}). It is far more challenging to construct a random ensemble that, on the one hand, satisfies a version of Theorem \ref{thm:JL}, and on the other is based on ``few random bits" or is constructed using a heavy-tailed random vector.

\vskip0.3cm
A significant breakthrough towards more general ``Johnson-Lindenstrauss transforms" came in \cite{MR2821584}, where it was shown that a matrix that satisfies a suitable version of the \emph{restricted isometry property}, can be converted to the wanted random ensemble by multiplying its columns by random signs. More accurately, let $\eps_1,...,\eps_n$ be independent, symmetric $\{-1,1\}$-valued random variables. Set $D_\eps={\rm diag}(\eps_1,...,\eps_n)$ and for a matrix $A : \R^n \to \R^m$ define
$$
A_\eps = AD_\eps.
$$
From here on we denote by $\Sigma_k$ the subset of $S^{n-1}$ consisting of vectors that are supported on at most $k$ coordinates.
\begin{Definition} \label{def:RIP1}
A matrix $A:\R^n \to \R^m$ satisfies the restricted isometry property of order $k$ and level $\delta \in (0,1)$ if
$$
\sup_{x \in \Sigma_k} \left| \|Ax\|_2^2 - 1 \right| \leq \delta.
$$
\end{Definition}

\begin{Theorem} \cite{MR2821584} \label{thm:KW}
There are absolute constants $c_0$ and $c_1$ such that the following holds.
Let $\lambda>0$ and $\rho \in (0,1)$. Consider $T \subset \R^n$ and let $k \geq c\log(e|T|/\lambda)$. If $A$ satisfies the restricted isometry property of order $k$ and at level $\delta<\rho/4$, then with probability at least $1-\lambda$, for every $t \in T$,
$$
(1-\rho)\|t\|_2^2 \leq  \|A_\eps t\|_2^2 \leq (1-\rho)\|t\|_2^2.
$$
\end{Theorem}

While Theorem \ref{thm:KW} does not recover the probability estimate from Theorem \ref{thm:JL}, it does imply at the constant probability level that $A_\eps$ is an almost isometry in the random ensemble sense: if $A$ is a matrix that $1 \pm \delta$-preserves the norms of vectors that are $c\log |T|$ sparse, then a typical realization of the random ensemble $A_\eps$, $1 \pm c^\prime \delta$ preserves the norms of all the elements in $T$.

\vskip0.3cm

Various extensions of Theorem \ref{thm:JL} that hold for arbitrary subsets of $\R^n$ have been studied over the years. In such extensions the ``complexity parameter" $\log |T|$ is replaced by more suitable counterparts. A rather general version of Theorem \ref{thm:JL} follows from a functional Bernstein inequality (see, e.g., \cite{MR3354613,MR3565471,Bed}), and to formulate that inequality in the gaussian case we require the following definition.
\begin{Definition}
Let $g_1,...g_n$ be independent, standard gaussian random variables. For $T \subset \R^n$ set
$$
\ell_*(T) = \E \sup_{t \in T} \left|\sum_{i=1}^n g_i t_i \right| \ \ \ {\rm and} \ \ \
d_T = \sup_{t \in T} \|t\|_2.
$$
Let
$$
\left(\frac{\ell_*(T)}{d_T}\right)^2
$$
be the \emph{critical dimension} of the set $T$.
\end{Definition}
The critical dimension appears naturally when studying the geometry of convex sets---for example, in the context of the Dvoretzky-Milman Theorem (see \cite{MR3331351} and references therein for more details). It is the natural alternative to $\log |T|$---which was suitable for finite subsets of sphere $S^{n-1}$.

\vskip0.3cm

Let $G=(g_i)_{i=1}^n$ be the standard gaussian random vector in $\R^n$, set $G_1,...,G_m$ to be independent copies of $G$ and put
$$
\Gamma=\sum_{i=1}^m \inr{G_i,\cdot}e_i
$$
to be the random ensemble used in Theorem \ref{thm:JL}.
\begin{Theorem} \label{thm:func-bern}
There exist absolute constants $c_0,c_1$ and $C$ such that the following holds.  If $T \subset \R^n$ and $u \geq c_0$ then with probability at least
$$
1-2\exp\left(-c_1u^2 \left(\frac{\ell_*(T)}{d_T}\right)^2 \right),
$$
for every $t \in T$,
\begin{equation} \label{eq:Bern-gaussian}
\left| \left\|m^{-1/2}\Gamma t\right\|_2^2 - \|t\|_2^2 \right| \leq C \left(ud_T \frac{\ell_*(T)}{\sqrt{m}} + u^2\left(\frac{\ell_*(T)}{\sqrt{m}}\right)^2\right).
\end{equation}
\end{Theorem}
One may use Theorem \ref{thm:func-bern} to ensure that the uniform error in \eqref{eq:Bern-gaussian} is at most $\max\{\rho,\rho^2\} d_T^2$. Indeed, if
$$
\frac{\ell_*(T)/d_T}{\sqrt{m}} \sim \rho,
$$
then with probability at least $1-2\exp(-c_3 \rho^2 m)$,
\begin{equation} \label{eq:JL-via-Bern}
\sup_{t \in T} \left| \left\|m^{-1/2}\Gamma t \right\|_2^2 - \|t\|_2^2 \right| \leq \max\{\rho,\rho^2\} d_T^2,
\end{equation}
which is a natural counterpart of Theorem \ref{thm:JL} once $\log |T|$ is replaced by $(\ell_*(T)/d_T)^2$.

\vskip0.3cm

As it happens, a version of Theorem \ref{thm:KW} that is analogous to \eqref{eq:JL-via-Bern} was proved in \cite{MR4023768}, using the notion of a multi-level RIP.
\begin{Definition} \label{def:Multi-RIP-old}
Let $L=\lceil \log_2 n \rceil$. For $\delta>0$ and $s \geq 1$ the matrix $A$ satisfies a multi-scale RIP with distortion $\delta$ and sparsity $s$ if, for every $1 \leq \ell \leq L$ and every $x \in \Sigma_{2^\ell s}$, one has
$$
\left| \|Ax\|_2^2 - \|x\|_2^2 \right| \leq \max\left\{2^{\ell/2}\delta, 2^\ell \delta^2 \right\}.
$$
\end{Definition}

Definition \ref{def:Multi-RIP-old} implies that if $k \geq s$ then
$$
\sup_{x \in \Sigma_k} \left|\|Ax\|_2^2-\|x\|_2^2\right| \leq \max\{ \sqrt{k} \delta, k\delta^2\}.
$$

\begin{Example} \label{Ex:multi-scale-gaussian}
Let $\Gamma:\R^n \to \R^m$ be a gaussian matrix as above and set $A=m^{-1/2}\Gamma$. It is standard to verify (using, for example, Theorem \ref{thm:func-bern} and a well-known estimate on $\ell_*(\Sigma_k)$) that with probability at least $1-2\exp(-ck\log(en/k))$,
$$
\sup_{x \in \Sigma_k} \left| \|Ax\|_2^2 - 1 \right| \leq C\sqrt{\frac{k\log(en/k)}{m}}.
$$
By the union bound over $k$ it follows that with a nontrivial probability, $A$ satisfies a multi-scale RIP with $s=1$ and $\delta \sim m^{-1/2}\sqrt{\log(en)}$. Observe that the second term in the multi-scale RIP---, namely $k\delta^2$, is not needed here.
\end{Example}

\begin{Remark}
Example \ref{Ex:multi-scale-gaussian} gives a good intuition on the role $\delta$ has in well-behaved situations: it should scale (roughly) like $1/\sqrt{m}$, where $m$ is the number of rows of the matrix $A$.
\end{Remark}

\vskip0.3cm
The following theorem is the starting point of this note: an estimate on the error a typical realization of the random ensemble $A_\eps=AD_\eps$ has when acting on an arbitrary $T \subset \R^n$, given that $A$ satisfies an appropriate multi-scale RIP.

\begin{Theorem} \cite{MR4023768} \label{thm:ORS}
There are absolute constants $c$ and $C$ such that the following holds.
Let $\eta,\rho>0$ and $A:\R^n \to \R^m$ that satisfies a multi-scale RIP with sparsity level $s=c(1+\eta)$ and distortion
\begin{equation} \label{eq:delta-thm-ORS}
\delta= C\frac{\rho d_T}{\max\{\ell_*(T),d_T\}}.
\end{equation}
Then for $T \subset \R^n$, with probability at least $1-\eta$,
\begin{equation} \label{eq:delta-thm-ORS-1}
\sup_{t \in T} \left| \|A_\eps t\|_2^2 - \|t\|_2^2 \right| \leq \max\{\rho^2,\rho\} d_T^2.
\end{equation}
\end{Theorem}

To put Theorem \ref{thm:ORS} in some context, if the belief is that $A_\eps$ should exhibit the same behaviour as the gaussian matrix $m^{-1/2}\Gamma$, then (keeping in mind that $\delta$ should scale like $1/\sqrt{m}$), ``a gaussian behaviour" as in Theorem \ref{thm:func-bern} is that with high probability,
$$
\sup_{t \in T} \left| \|A_\eps t\|_2^2 - \|t\|_2^2 \right| \leq c\left(\delta d_T \ell_*(T) + (\delta \ell_*(T))^2\right).
$$
Observe that $\ell_*(T) \gtrsim d_T$, implying by \eqref{eq:delta-thm-ORS} that $\rho \sim \delta \ell_*(T)/d_T$. Hence, the error in \eqref{eq:delta-thm-ORS-1} in terms of $\delta$ is indeed
$$
\sim d_T \delta  \ell_*(T) + (\delta \ell_*(T))^2.
$$
However, despite the ``gaussian error", the probability estimate in Theorem \ref{thm:ORS} is far weaker than in Theorem \ref{thm:func-bern}---it is just at the constant level.

\vskip0.3cm

Our main result is that using a modified, seemingly less restrictive version of the multi-scale RIP, $A_\eps$ acts on $T$ as if it were a gaussian operator: achieving the same distortion and probability estimate as in Theorem \ref{thm:func-bern}.

\begin{Definition}
Let $A:\R^n \to \R^m$ be a matrix. For $\delta>0$ let $1 \leq k^* \leq n$ be the largest such that for every $1 \leq k \leq k_*$, $A$ is a $(\delta,k)$ regular; that is, for every $1 \leq k \leq k_*$
$$
\sup_{x \in \Sigma_k} \left| \|Ax\|_2^2 - \|x\|_2^2 \right| \leq \delta \sqrt{k}.
$$
\end{Definition}

\begin{Theorem} \label{thm:main-intro}
There exist absolute constants $c$ and $C$ such that the following holds. Let $\delta>0$ and set $\Lambda=\max\{1,\delta^2 \log(n\delta^2)\}$. If $k_* \geq 1/\delta^2$, $T \subset \R^n$ and $u \geq 1$, then with probability at least
$$
1-2\exp\left(-cu^2 \left(\frac{\ell_*(T)}{d_T} \right)^2 \right),
$$
we have that
\begin{equation} \label{eq:main-error}
\sup_{t \in T} \left| \|A_\eps t \|_2^2 - \|t\|_2^2 \right| \leq Cu^2 \left( \Lambda \cdot d_T \delta \ell_*(T)+(\delta \ell_*(T))^2\right).
\end{equation}
\end{Theorem}

\begin{Remark}
The sub-optimality in Theorem \ref{thm:main-intro} lies in the factor $\Lambda$---in the range where $\delta$ is relatively large: at least $1/\sqrt{\log n}$. For $\delta \leq 1/\sqrt{\log n}$ we have that $\Lambda=1$ and Theorem \ref{thm:main-intro} recovers the functional Bernstein inequality for $u \sim 1$; that holds despite the fact that $A_\eps$ is based only on $n$ ``random bits".

Moreover, for the error in \eqref{eq:main-error} to have a chance of being a nontrivial two-sided estimate, i.e., that for some $0<\rho<1$ and every $t \in T$,
$$
\left| \|A_\eps t\|_2^2 - \|t\|_2^2 \right| \leq \rho d_T^2,
$$
$\delta$ has to be smaller than $\sim d_T/\ell_*(T)$. In particular, if the critical dimension of $T$, $(\ell_*(T)/d_T)^2$ is at least $\log n$, a choice of $\delta \leq (d_T/\ell_*(T))$ leads to $\Lambda=1$ and thus to an optimal outcome in Theorem \ref{thm:main-intro}.
\end{Remark}
\vskip0.3cm
Theorem \ref{thm:main-intro} clearly improves the probability estimate from Theorem \ref{thm:ORS}. The other (virtual) improvement is that the matrix $A$ need only be $(\delta,k)$-regular for $k \leq 1/\delta^2$, and the way $A$ acts on $\Sigma_k$ for $k > 1/\delta^2$ is of no importance. The reason for calling that improvement ``virtual" is the following observation:
\begin{Lemma} \label{lemma:low-to-high-sparsity}
If $k^* \geq 1/\delta^2$ then for any $1 \leq s \leq n$,
$$
\sup_{x \in \Sigma_s} \left| \|Ax\|_2^2 - 1 \right| \leq 4\max\{\delta \sqrt{s}, \delta^2 s\}.
$$
\end{Lemma}
In other words, the second term in the multi-scale RIP condition follows automatically from the first one and the fact that $k^*$ is sufficiently large.
\vskip0.3cm

\proof Let $x \in \Sigma_s$ for $s \geq k_*$, and let $(J_i)_{i=1}^\ell$ be a decomposition of the support of $x$ to coordinate blocks of cardinality $k_*/2$. Set $y_i = P_{J_i} x$, that is, the projection of $x$ onto ${\rm span}\{e_m : m \in J_i\}$ and write $x = \sum_{i=1}^\ell y_i$. Note that $\ell \leq 4s/k_*$ and that
$$
\|A x\|_2^2 = \bigl\|\sum_{i=1}^\ell A y_i\bigr\|^2=\sum_{i=1}^\ell \|Ay_i\|_2^2 + \sum_{i \not = j} \inr{Ay_i,Ay_j}.
$$
The vectors $y_i$ are orthogonal and so $\|x\|_2^2 = \sum_{i=1}^\ell \|y_i\|_2^2$. Therefore,
$$
\left| \|Ax\|_2^2 - \|x\|_2^2 \right| \leq \sum_{i=1}^\ell \left| \|Ay_i\|_2^2 - \|y_i\|_2^2 \right| + \left|\sum_{i \not = j} \inr{Ay_i,Ay_j}\right|.
$$
For the first term, as each $y_i$ is supported on at most $k_*/2$ coordinates, it follows from the regularity condition that
$$
\sum_{i=1}^\ell \left| \|Ay_i\|_2^2 - \|y_i\|_2^2 \right| \leq \delta \sqrt{k_*/2} \sum_{i=1}^\ell \|y_i\|_2^2 = \delta \sqrt{k_*/2} \|x\|_2^2 \leq \delta \sqrt{s}.
$$
As for the second term, since $y_i$ and $y_j$ are orthogonal, $\|y_i+y_j\|_2 = \|y_i-y_j\|_2$ and
\begin{align*}
\inr{Ay_i,Ay_j} = &\frac{1}{4} \left(\|A(y_i+y_j)\|_2^2-\|A(y_i-y_j)\|_2^2\right)
\\
= & \frac{1}{4} \left(\|A(y_i+y_j)\|_2^2-\|y_i+y_j\|_2^2\right) - \frac{1}{4} \left(\|y_i-y_j\|_2^2 -\|A(y_i-y_j)\|_2^2\right).
\end{align*}
Thus, by the regularity of $A$ and as $|{\rm supp}(y_i \pm y_j)| \leq k_*$,
\begin{align*}
|\inr{Ay_i,Ay_j}| \leq & \frac{1}{4} \left(\delta \sqrt{k_*} \|y_i+y_j\|_2^2 + \delta \sqrt{k_*} \|y_i-y_j\|_2^2\right)
\\
\leq & \frac{1}{2} \delta \sqrt{k_*} (\|y_i\|_2^2 + \|y_j\|_2^2).
\end{align*}

Taking the sum over all pairs $i \not = j$, $i,j \leq \ell$, each factor $\|y_i\|_2^2$ appears at most $2\ell$ times, and $2\ell \leq 8s/k_*$. Hence, using that $1/\delta^2 \leq k_*$
$$
\sum_{i \not = j} |\inr{Ay_i,Ay_j}| \leq \frac{1}{2}\delta \sqrt{k_*} \cdot \frac{8s}{k_*} \sum_{i=1}^\ell \|y_i\|_2^2 \leq 2\delta \frac{s}{\sqrt{k_*}}\|x\|_2^2 \leq 4\delta^2 s.
$$
\endproof

\vskip0.3cm

Clearly, Theorem \ref{thm:main-intro} implies a suitable version of Theorem \ref{thm:ORS}.
\begin{Corollary} \label{cor:main-intro}
There exist absolute constants $c$ and $c_1$ such that the following holds. Let $A$ be as above, set $T \subset \R^n$ and $0<\delta<1/\sqrt{\log n}$. Let $\rho=c\delta \ell_*(T)/d_T$. Then with probability at least $1-2\exp(-c_1\rho^2/\delta^2)$,
$$
\sup_{t \in T} \left| \|A_\eps t \|_2^2 - \|t\|_2^2 \right| \leq \rho d_T^2.
$$
\end{Corollary}

\begin{Remark}
Recalling the intuition that $m \sim 1/\delta^2$, the outcome of Corollary \ref{cor:main-intro} coincides with the estimate in \eqref{eq:JL-via-Bern}.
\end{Remark}

\vskip0.3cm

In Section \ref{sec:appl} we present one simple application of Theorem \ref{thm:main-intro}. We show that column randomization of a typical realization of a Bernoulli circulant matrix (complete or partial) exhibits an almost gaussian behaviour (conditioned on the generating vector). In particular, only $2n$ random bits ($n$ from the generating Bernoulli vector and $n$ from the column randomization) are required if one wishes to create a random ensemble that is, effectively, an almost isometry.

\vskip0.3cm
The proof of Theorem \ref{thm:main-intro} is based on a chaining argument. For more information on the \emph{generic chaining} mechanism, see Talagrand's treasured manuscript \cite{MR3184689}. We only require relatively basic notions from generic chaining theory, as well as the celebrated \emph{majorizing measures theorem}.
\begin{Definition}
Let $T \subset \R^n$. A collection of subsets of $T$, $(T_s)_{s \geq 0}$, is an \emph{admissible sequence} if $|T_0|=1$ and for $s \geq 1$, $|T_s| \leq 2^{2^s}$. For every $t \in T$ denote by $\pi_st$ a nearest point to $t$ in $T_s$ with respect to the Euclidean distance. Set $\Delta_st = \pi_{s+1}t-\pi_st$ for $s \geq 1$ and let $\Delta_0 t = \pi_0t$.

The $\gamma_2$ functional with respect to the $\ell_2$ metric is defined by
$$
\gamma_2(T,\| \ \|_2) = \inf_{(T_s)} \sup_{t \in T} \sum_{s \geq 0} 2^{s/2}\|\Delta_s t\|_2,
$$
where the infimum is taken with respect to all admissible sequences of $T$.
\end{Definition}

An application of Talagrand's \emph{majorizing measures theorem} to the gaussian process $t \to \sum_{i=1}^n g_i t_i$ shows that $\gamma_2(T,\| \ \|_2)$ and $\ell_*(T)$ are equivalent:
\begin{Theorem} \label{thm:MM}
There are absolute constants $c$ and $C$ such that for every $T \subset \R^n$,
$$
c\gamma_2(T,\| \ \|_2) \leq \ell_*(T) \leq C\gamma_2(T,\| \ \|_2).
$$
\end{Theorem}
The proof of Theorem \ref{thm:MM} can be found, for example, in \cite{MR3184689}.

\section{Proof of Theorem \ref{thm:main-intro}}
We begin the proof with a word about notation: throughout, absolute constants, that is, positive numbers that are independent of all the parameters involved in the problem, are denoted by $c$, $c_1$, $C$, etc. Their value may change from line to line.

\vskip0.3cm

As noted previously, the proof is based on a chaining argument. Let $(T_s)_{s \geq 0}$ be an optimal admissible sequence of $T$. Set $s_0$ to satisfy that $2^{s_0}$ is the critical dimension of $T$, i.e.,
$$
2^{s_0} = \left(\frac{\ell_*(T)}{d_T}\right)^2
$$
(without loss of generality we may assume that equality holds). Let $s_0 \leq s_1$ to be named in what follows and observe that
$$
\|A_\eps t\|_2^2 = \|A_\eps (t-\pi_{s_1}t) + A_\eps \pi_{s_1}t \|_2^2 = \|A_\eps (t-\pi_{s_1}t)\|_2^2 + 2\inr{A_\eps (t-\pi_{s_1}t),A_\eps \pi_{s_1}t} + \|A_\eps \pi_{s_1} t\|_2^2.
$$
Writing $t-\pi_{s_1}t=\sum_{s \geq s_1} \Delta_s t$, it follows that
\begin{equation} \label{eq:est1}
\|A_\eps (t-\pi_{s_1}t)\|_2  \leq \sum_{s \geq s_1} \|A_\eps \Delta_s t\|_2,
\end{equation}
and
$$
\left|\inr{A_\eps (t-\pi_{s_1}t),A_\eps \pi_{s_1}t}\right| \leq \|A_\eps (t-\pi_{s_1}t)\|_2 \cdot \|A_\eps \pi_{s_1}t\|_2;
$$
Therefore, setting
$$
\Psi^2=\sup_{t \in T} \left| \|A_\eps \pi_{s_1} t\|_2^2 - \|\pi_{s_1}t\|_2^2 \right| \ \ \ {\rm and} \ \ \ \Phi=\sum_{s \geq s_1} \|A_\eps \Delta_s t\|_2
$$
we have that for every $t \in T$,
$$
\|A_\eps \pi_{s_1} t\|_2 \leq \sqrt{\Psi^2+d_T^2}
$$
and
\begin{equation} \label{eq:est2}
\left|\inr{A_\eps (t-\pi_{s_1}t),A_\eps \pi_{s_1}t}\right| \leq \Phi \cdot \sqrt{\Psi^2+d_T^2}.
\end{equation}

Hence,
\begin{equation} \label{eq:est3}
\sup_{t \in T} \left| \|A_\eps t\|_2^2 - \|t\|_2^2 \right| \leq \Psi^2 + 2 \Phi \cdot \sqrt{\Psi^2+d_T^2}  + \Phi^2+\sup_{t \in T} \left| \|\pi_{s_1}t\|_2^2 - \|t\|_2^2 \right|.
\end{equation}

To estimate the final term, note that for every $t \in T$,
\begin{align*}
\left|\|t\|_2^2 - \|\pi_{s_1}t\|_2^2 \right| \leq & \|t-\pi_{s_1}t\|_2^2 + 2 \left|\inr{t-\pi_{s_1}t,\pi_{s_1}t}\right| \leq \|t-\pi_{s_1}t\|_2^2 + 2  \|t-\pi_{s_1}t\|_2 \cdot \|\pi_{s_1}t\|_2
\\
\leq & \bigl(\sum_{s \geq s_1} \|\Delta_s t\|_2\bigr)^2 + 2d_T \sum_{s \geq s_1} \|\Delta_s t\|_2.
\end{align*}
By the definition of the $\gamma_2$ functional and the majorizing measures theorem, for every integer $s$ and every $t \in T$,
$$
\|\Delta_st\|_2 \leq 2^{-s/2}\gamma_2(T,\| \ \|_2) \leq c_12^{-s/2}\ell_*(T).
$$
Thus,
$$
\sum_{s \geq s_1} \|\Delta_s t\|_2 \leq c_1\ell_*(T) \sum_{s \geq s_1} 2^{-s/2} \leq c_2 2^{-s_1/2} \ell_*(T),
$$
and
\begin{equation} \label{eq:est4}
\left|\|\pi_{s_1}t\|_2^2 - \|t\|_2^2\right| \leq c_3 \left(2^{-s_1} \ell_*^2(T) + d_T 2^{-s_1/2}\ell_*(T) \right).
\end{equation}

Equation \eqref{eq:est4} and the wanted estimate in Theorem \ref{thm:main-intro} hint on the identity of $2^{s_1}$: it should be larger than $1/\delta^2$. Recalling that $s_0 \leq s_1$ and that $2^{s_0}=(\ell_*(T)/d_T)^2$, set
$$
2^{s_1} = \max\left\{\frac{1}{\delta^2}, \left(\frac{\ell_*(T)}{d_T}\right)^2, \log\left(e(1+n\delta^2)\right) \right\}.
$$
The reason behind the choice of the third term will become clear in what follows.

With that choice of $s_1$,
\begin{equation} \label{eq:est5}
\sup_{t \in T} \left|\|\pi_{s_1}t\|_2^2 - \|t\|_2^2\right| \leq c_3 \left( \delta^2 \ell_*^2(T) + d_T \delta \ell_*(T) \right),
\end{equation}
and the nontrivial part of the proof is to control $\Phi$ and $\Psi$  with high probability that would lead to the wanted estimate on \eqref{eq:est3}.

\subsection{A decoupling argument}
For every $t \in \R^n$,
\begin{equation} \label{eq:expand}
\|A_\eps t\|_2^2 = \sum_{i,j} \inr{Ae_i,Ae_j} \eps_i \eps_j t_i t_j = \sum_{i=1}^n \|Ae_i\|_2^2 t_i^2 + \sum_{i \not =j } \inr{Ae_i,Ae_j} \eps_i \eps_j t_i t_j.
\end{equation}
By the assumption that $A$ is $(\delta,1)$-regular,
$$
\max_{1 \leq i \leq n} \left|\|A e_i\|_2^2 - 1\right| \leq \delta,
$$
and noting that $d_T \leq c \ell_*(T)$,
$$
\left|\sum_{i=1}^n \|Ae_i\|_2^2 t_i^2 - \|t\|_2^2 \right| = \left|\sum_{i=1}^n \left(\|Ae_i\|_2^2-1\right) t_i^2\right| \leq \delta \|t\|_2^2 \leq \delta d_T^2 \leq cd_T \delta \ell_*(T).
$$
Next, we turn to the ``off-diagonal" term in \eqref{eq:expand}. For $t \in \R^n$ let
$$
Z_t = \sum_{i \not = j} \eps_i \eps_j t_i t_j \inr{Ae_i,Ae_j}
$$
and let us obtain high probability estimates on $\sup_{u \in U} |Z_u|$ for various sets $U$.
\vskip0.3cm
The first step in that direction is decoupling: let $\eta_1,...,\eta_n$ be independent $\{0,1\}$-valued random variables with mean $1/2$. Set $I=\{i : \eta_i =1\}$ and observe that for every $(\eps_i)_{i=1}^n$,
\begin{equation} \label{eq:decoupling-est}
\sup_{u \in U} |Z_u| \leq  4  \sup_{u \in U} \E_\eta \left|\inr{A \left(\sum_{i \in I} \eps_i u_i e_i\right), A \left(\sum_{j \in I^c} \eps_j u_j e_j\right)}\right|.
\end{equation}
Indeed, for every $(\eps_i)_{i=1}^n$,
\begin{align*}
\sup_{u \in U} \left|\sum_{i \not = j} \inr{Ae_i,Ae_j}\eps_i \eps_j u_i u_j \right| = & 4 \sup_{u \in U} \left| \sum_{i \not = j} \inr{Ae_i,Ae_j}\E_\eta (1-\eta_i) \eta_j  \cdot \eps_i\eps_j u_i u_i \right|
\\
\leq & 4  \sup_{u \in U} \E_\eta \left| \sum_{i \in I, \ j \in I^c} \inr{Ae_i,Ae_j}\eps_i\eps_j u_i u_j \right|;
\end{align*}
and for every $u \in \R^n$,
$$
\sum_{i \in I, \ j \in I^c} \inr{Ae_i,Ae_j}\eps_i\eps_j u_i u_j  = \inr{A \bigl(\sum_{j \in I^c} \eps_ju_j e_j\bigr), A\bigl(\sum_{i \in I}  \eps_i u_i e_i \bigr)}.
$$
Equation \eqref{eq:decoupling-est} naturally leads to the following definition:
\begin{Definition}
For $v \in \R^n$ and $I \subset \{1,...,n\}$, set
\begin{equation} \label{eq:W-v}
W_{v,I}=A^* A \bigl(\sum_{i \in I} \eps_i v_i e_i \bigr).
\end{equation}
\end{Definition}

Recall that $\pi_{s+1}t$ is the nearest point to $t$ in $T_{s+1}$ and $\Delta_s t =\pi_{s+1}t-\pi_st$.
\begin{Lemma} \label{lemma:key-observe-chaining}
For every $t$ and every $(\eps_i)_{i=1}^n$,
\begin{align} \label{eq:chaining-1}
|Z_{\pi_{s_1}t}| \leq & 4 \sum_{s =s_0}^{s_1-1} \left( \E_\eta \left| \sum_{j \in I^c} \eps_j (\Delta_s t)_j (W_{\pi_{s+1} t,I})_j\right|+\E_\eta \left|\sum_{i \in I} \eps_i (\Delta_s t)_i (W_{\pi_{s} t,I^c})_i\right| \right) \nonumber
\\
+ & 4 \E_\eta \left|  \sum_{j \in I^c} \eps_j (\pi_{s_0} t)_j (W_{\pi_{s_0} t,I})_j\right|.
\end{align}
\end{Lemma}

\proof Fix an integer $s$. With the decoupling argument in mind, fix $I \subset \{1,...,n\}$ and observe that
\begin{align*}
& \inr{A \bigl(\sum_{i \in I} \eps_i (\pi_{s+1} t)_i e_i\bigr), A \bigl(\sum_{j \in I^c} \eps_j (\pi_{s+1} t)_j e_j\bigr)}
\\
= & \inr{A \bigl(\sum_{i \in I} \eps_i (\pi_{s+1} t)_i e_i\bigr), A \bigl(\sum_{j \in I^c} \eps_j (\Delta_s t)_j e_j\bigr)}+\inr{A \bigl(\sum_{i \in I} \eps_i (\pi_{s+1} t)_i e_i\bigr), A \bigl(\sum_{j \in I^c} \eps_j (\pi_{s} t)_j e_j\bigr)}
\\
= & \inr{A \bigl(\sum_{i \in I} \eps_i (\pi_{s+1} t)_i e_i\bigr), A \bigl(\sum_{j \in I^c} \eps_j (\Delta_s t)_j e_j\bigr)}+\inr{A \bigl(\sum_{i \in I} \eps_i (\Delta_{s} t)_i e_i\bigr), A \bigl(\sum_{j \in I^c} \eps_j (\pi_{s} t)_j e_j\bigr)}
\\
+ & \inr{A \bigl(\sum_{i \in I} \eps_i (\pi_{s} t)_i e_i\bigr), A \bigl(\sum_{j \in I^c} \eps_j (\pi_{s} t)_j e_j\bigr)}.
\end{align*}
Moreover,
\begin{align*}
\inr{A \bigl(\sum_{i \in I} \eps_i (\pi_{s+1} t)_i e_i\bigr), A \bigl(\sum_{j \in I^c} \eps_j (\Delta_s t)_j e_j\bigr)} = & \inr{ A^*A \bigl(\sum_{i \in I} \eps_i (\pi_{s+1} t)_i e_i\bigr), \bigl(\sum_{j \in I^c} \eps_j (\Delta_s t)_j e_j\bigr)}
\\
= & \sum_{j \in I^c} \eps_j (\Delta_s t)_j (W_{\pi_{s+1} t,I})_j.
\end{align*}
and
\begin{equation*}
\inr{A \bigl(\sum_{i \in I} \eps_i (\Delta_{s} t)_i e_i\bigr), A \bigl(\sum_{j \in I^c} \eps_j (\pi_{s} t)_j e_j\bigr)}=\sum_{i \in I} \eps_i (\Delta_s t)_i (W_{\pi_{s} t,I^c})_i.
\end{equation*}

Combining these observations, for every $t \in T$ and $(\eps_i)_{i=1}^n$,
\begin{align*}
\frac{1}{4} Z_{\pi_{s_1}t} = & \E_\eta \inr{A \bigl(\sum_{i \in I} \eps_i (\pi_{s_1}t)_i e_i \bigr), A \bigl(\sum_{j \in I^c} \eps_j (\pi_{s_1}t)_j e_j\bigr)}
\\
= & \E_\eta \sum_{s=s_0}^{s_1-1} \sum_{i \in I} \eps_i (\Delta_s t)_i (W_{\pi_s t,I^c})_i
\\
+ & \E_\eta \sum_{s=s_0}^{s_1-1} \sum_{j \in I^c} \eps_j (\Delta_s t)_j (W_{\pi_{s+1} t,I})_j
\\
+ & \E_\eta \sum_{j \in I^c} \eps_j (\pi_{s_0} t)_j (W_{\pi_{s_0} t,I})_j,
\end{align*}
from which the claim follows immediately.
\endproof

As part of the decoupling argument and to deal with the introduction of the random variables $(\eta_i)_{i=1}^n$ in \eqref{eq:decoupling-est}, we will use the following elementary fact:
\begin{Lemma} \label{rem:simple-Fubini}
Let $f$ be a function of $(\eps_i)_{i=1}^n$ and $(\eta_i)_{i=1}^n$. If $\E_\eta \E_\eps  |f|^q \leq \kappa^q$ then with $(\eps_i)_{i=1}^n$-probability at least $1-\exp(-q)$,
\begin{equation} \label{eq:simple-Fubini}
\E_\eta (|f| \ | (\eps_i)_{i=1}^n) \leq e\kappa.
\end{equation}
\end{Lemma}

\proof For a nonnegative function $h$ we have that point-wise, $\IND_{\{h \geq t\}} \leq h^q/t^q$. Let $h=\E_\eta (|f| \ | (\eps_i)_{i=1}^n)$ and note that
$$
Pr_\eps \left(\E_\eta (|f| \ | (\eps_i)_{i=1}^n) \geq u\kappa \right) =\E_\eps \IND_{\{h \geq u\kappa\}} \leq (u\kappa)^{-q}(\E_\eps h^q).
$$
By Jensen's inequality followed by Fubini's Theorem,
$$
\E_\eps h^q = \E_\eps (\E_\eta (|f| \ | (\eps_i)_{i=1}^n)^q \leq \E_\eps \E_\eta |f|^q \leq \kappa^q,
$$
and setting $u=e$ proves \eqref{eq:simple-Fubini}.
\endproof
\vskip0.3cm
We shall use Lemma \ref{rem:simple-Fubini} in situations where we actually have more information---namely that for {\emph any} $(\eta_i)_{i=1}^n$, $\E_\eps (|f|^q \ |(\eta_i)_{i=1}^n) \leq \kappa^q$ for a well chosen $\kappa$. As a result,
$$
\E_\eta \bigl( |f| \big| (\eps_i)_{i=1}^n \bigr) \leq e \kappa \ \ {\rm with \ probability \ at \ least \ } 1-\exp(-q).
$$
\vskip0.3cm

Taking into account Lemma \ref{lemma:key-observe-chaining} and Lemma \ref{rem:simple-Fubini}, it follows that if one wishes to estimate $\sup_{t \in T} |Z_{\pi_{s_1}t}|$ using a chaining argument, it suffices to obtain, for every $I \subset \{1,...,n\}$, bounds on moments of random variables of the form
\begin{equation} \label{eq:needed-rv}
\sum_{j \in I^c} \eps_j (\Delta_s t)_j (W_{\pi_{s+1} t,I})_j, \ \ \  \sum_{i \in I} \eps_i (\Delta_s t)_i (W_{\pi_{s} t,I^c})_i, \ \ \ {\rm and} \ \ \  \sum_{j \in I^c} \eps_j (\pi_{s_0} t)_j (W_{\pi_{s_0} t,I})_j,
\end{equation}
as that results in high probability estimates on each of the terms in \eqref{eq:chaining-1}. And as there are at most $2^{2^{s+3}}$ random variables involved in this chaining argument at the $s$-stage, the required moment in \eqref{eq:needed-rv} is $q \sim 2^s$ for the first two terms and $q \sim 2^{s_0}$ for the third one.

\subsection{Preliminary estimates}
For $J \subset \{1,...,n\}$ let $P_Jx = \sum_{j \in J} x_je_j$ be the projection of $x$ onto ${\rm span}(e_j)_{j \in J}$. The key lemma in the proof of Theorem \ref{thm:main-intro} is:

\begin{Lemma} \label{lemma:key-observation-1}
There exists an absolute constant $c$ such that the following holds.
Let $I \subset \{1,...,n\}$ and $W_{v,I}$ be as in \eqref{eq:W-v}. Set $J \subset I^c$ such that $|J| \leq k_*$. Then for $q \geq |J|$,
$$
\left(\E \|P_J W_{v,I}\|_2^q\right)^{1/q} \leq c\sqrt{q}\delta \|v\|_2.
$$
\end{Lemma}

\proof Let $S^J$ be the Euclidean unit sphere in the subspace ${\rm span}(e_j)_{j \in J}$ and let $U$ be a maximal $1/10$-separated subset of $S^J$. By a volumetric estimate (see, e.g. \cite{MR3331351}), there is an absolute constant $c_0$ such that $|U| \leq \exp(c_0 |J|)$. Moreover, a standard approximation argument shows that for every $y \in \R^n$,
$$
\|P_J y\|_2 = \sup_{x \in S^J} \inr{y,x} \leq c_1 \max_{x \in U} \inr{y,x},
$$
where $c_1$ is a suitable absolute constant. Therefore,
$$
\sup_{x \in S^J} \inr{W_{v,I},x} \leq c_1 \max_{x \in U} \inr{W_{v,I},x},
$$
and it suffices to control, with high probability,
$$
\max_{x \in U} \inr{\sum_{i \in I} \eps_i v_i e_i, A^*Ax} = \max_{x \in U} \sum_{i \in I} \eps_i v_i \left(A^*Ax\right)_i.
$$
Fix $x \in U$, recall that $A$ is $(\delta,k)$-regular for $1 \leq k \leq k_*$ and we first explore the case $1 \leq q \leq k_*/2$.

Denote by $I^\prime \subset I$ the set of indices corresponding to the $q$ largest values of $|(A^*Ax)_i|, \ i \in I$. If $|I| \leq q$ then set $I^\prime = I$.

It is straightforward to verify (e.g., using H\"{o}ffding's inequality) that there is an absolute constant $c_2$ such that
\begin{align*}
\left\|\sum_{i \in I} \eps_i v_i \left(A^*Ax\right)_i\right\|_{L_q} \leq & \sum_{i \in I^\prime} |v_i \left(A^*Ax\right)_i| + c_2\sqrt{q} \bigl(\sum_{i \in I \backslash I^\prime} v_i^2 \left(A^*Ax\right)_i^2 \bigr)^{1/2}
\\
\leq & \|v\|_2 \|P_{I^\prime} (A^*Ax)\|_2 + c_2\sqrt{q} \cdot \frac{\|P_{I^\prime} (A^*Ax)\|_2}{\sqrt{q}}\|v\|_2
\\
\leq & c_3\|v\|_2 \|P_{I^\prime} (A^*Ax)\|_2,
\end{align*}
where we used that fact that for $i \in I \backslash I^\prime$,
$$
|\left(A^*Ax\right)_i| \leq \frac{\|P_{I^\prime} (A^*Ax)\|_2}{\sqrt{q}}.
$$
Therefore,
$$
\left\|\sum_{i \in I} \eps_i v_i \left(A^*Ax\right)_i\right\|_{L_q} \leq c_3\|v\|_2 \max_{I^\prime \subset I, \ |I^\prime|=q} \sup_{z \in S^{I^\prime}} \inr{A^*Ax,z}.
$$
Note that $x$ is supported in $J \subset I^c$, while each `legal' $z$ is supported in a subset of $I$; in particular, $x$ and $z$ are orthogonal, implying that for every such $z$,
$$
\|x+z\|_2=\|x-z\|_2 \leq 2 \ \ {\rm and} \ \ |{\rm supp}(x+z)|, |{\rm supp}(x-z)| \leq 2q.
$$
Thus, by the $(\delta,2q)$-regularity of $A$ (as $2q \leq k_*$),
\begin{align*}
\left|\inr{A^*Ax,z}\right| = & \left|\inr{Az,Ax}\right| = \frac{1}{4}\left| \|A(x+z)\|_2^2 - \|A(x-z)\|_2^2\right|
\\
= & \frac{1}{4}\left| \left(\|A(x+z)\|_2^2-\|x+z\|_2^2\right) - \left(\|A(x-z)\|_2^2-\|x-z\|_2^2 \right) \right|
\\
\leq & \frac{\delta}{4} \sqrt{2q} \cdot \max\{\|x+z\|_2^2, \|x-z\|_2^2\} \leq \delta \sqrt{q},
\end{align*}
and it follows that
$$
\left\|\sum_{i \in I} \eps_i v_i \left(A^*Ax\right)_i\right\|_{L_q} \leq c_4\|v\|_2 \delta \sqrt{q}.
$$

Turning to the case $q \geq k_*/2$, let $I^\prime$ be the set of indices corresponding to the $k_*/2$ largest coordinates of $(|(A^*Ax)_i|)_{i \in I}$. Therefore,

\begin{align} \label{eq:large-q}
\left\| \sum_{i \in I} \eps_i v_i (A^*Ax)_i \right\|_{L_q} \leq &
\sum_{i \in I^\prime} v_i (A^*Ax)_i + c_2 \sqrt{q} \left(\sum_{i \in I \backslash I^\prime} v_i^2 (A^*Ax)_i^2\right)^{1/2}  \nonumber
\\
\leq & \|v\|_2 \cdot \|P_{I^\prime} A^*Ax\|_2 \left(1+c_5\sqrt{\frac{q}{k_*}}\right),
\end{align}
using that for $i \in I \backslash I^\prime$,
$$
|(A^*Ax)_i| \leq \frac{\|P_{I^\prime} A^*Ax\|_2}{\sqrt{k_*/2}}.
$$
Recall that $|{\rm supp}(x)|=|J| \leq k_*$ and that $|I^\prime|=k_*/2$. The same argument used previously shows that
$$
\|P_{I^\prime} A^*Ax\|_2 \leq c_6\delta \sqrt{k_*};
$$
hence,
$$
\left\|\sum_{i \in I} \eps_i v_i \left(A^*Ax\right)_i\right\|_{L_q} \leq c_{7}\|v\|_2 \delta \sqrt{q},
$$
and the estimate holds for each $q \geq |J|$ for that fixed $x$.  

Setting $u \geq 1$, it follows from Chebychev's inequality that with probability at least $1-2\exp(-c_8u^2q)$,
$$
\left|\sum_{i \in I} \eps_i v_i \left(A^*Ax\right)_i\right| \leq c_9  u \|v\|_2 \delta \sqrt{q},
$$
and by the union bound, recalling that $q \geq |J|$, the same estimate holds uniformly for every $x \in U$, provided that $u \geq c_{10}$. Thus, with probability at least $1-2\exp(-cu^2q)$,
$$
\|P_J W_{v,I}\|_2 \leq c^\prime \max_{x \in U} \left| \sum_{i \in I} \eps_i v_i (A^*Ax)_i \right| \leq Cu\|v\|_2 \delta \sqrt{q},
$$
and the wanted estimate follows from tail integration.
\endproof

The next observation deals with more refined estimates on random variables of the form
$$
X_{a,b}=\sum_{i \in I^c} \eps_i a_i b_i.
$$
Once again, we use the fact that for any $J \subset I^c$
\begin{equation} \label{eq:X-a-b}
\|X_{a,b}\|_{L_q} \leq \sum_{j \in J} |a_ib_i| + c\sqrt{q}\bigl(\sum_{j \in I^c \backslash J} a_i^2b_i^2\bigr)^{1/2}.
\end{equation}

As one might have guessed, the choice of $b$ will be vectors of the form $W_{\pi_s t ,I}$. These are random vectors that are independent of the Bernoulli random variables involved in the definition of $X$. At the same time, $a$ will be deterministic.

\vskip0.3cm

Without loss of generality assume that $a_i,b_i \geq 0$ for every $i$. Let $J_1$ be the set of indices corresponding to the $k_1$ largest coordinates of $(a_i)_{i \in I^c}$, $J_2$ is the set corresponding to the following $k_2$ largest coordinates, and so on. The choice of $k_1,k_2,...,$ will be specified in what follows.

Note that for any $\ell>1$,
$$
\|P_{J_\ell} a\|_\infty \leq \frac{\|P_{J_{\ell-1}}a\|_2}{\sqrt{|J_{\ell-1}|}} = \frac{\|P_{J_{\ell-1}}a\|_2}{\sqrt{k_{\ell-1}}}.
$$
Set $J=J_1$, and by \eqref{eq:X-a-b},
\begin{align} \label{eq:X-a-b-2}
\|X_{a,b}\|_{L_q} \leq & \sum_{j \in J_1} |a_ib_i| + c\sqrt{q}\bigl(\sum_{\ell \geq 1} \sum_{j \in J_{\ell+1}} a_i^2b_i^2\bigr)^{1/2} \nonumber
\\
\leq & \|P_{J_1} a\|_{2} \|P_{J_1} b\|_2 + c\sqrt{q} \bigl(\sum_{\ell \geq 1} \|P_{J_{\ell+1}} a\|_\infty^2 \|P_{J_{\ell+1}} b\|_2^2 \bigr)^{1/2} \nonumber
\\
\leq & \|P_{J_1} a\|_{2} \|P_{J_1} b\|_2 + c\sqrt{q} \bigl(\sum_{\ell \geq 1} \frac{\|P_{J_\ell} a\|_2^2}{|J_\ell|} \|P_{J_{\ell+1}} b\|_2^2 \bigr)^{1/2} \nonumber
\\
\leq & \|P_{J_1} a\|_{2} \|P_{J_1} b\|_2 +c \sqrt{q} \max_{\ell \geq 2} \frac{\|P_{J_\ell} b\|_2}{\sqrt{k_{\ell-1}}} \cdot \bigl(\sum_{\ell \geq 2} \|P_{J_\ell} a\|_2^2 \bigr)^{1/2} \nonumber
\\
\leq & \|a\|_{2} \Bigl(\|P_{J_1} b\|_2 + c\sqrt{q} \max_{\ell \geq 2} \frac{\|P_{J_{\ell}} b\|_2}{\sqrt{k_{\ell-1}}} \Bigr).
\end{align}

And, in the case where $|J_\ell|=k$ for every $\ell$, it follows that
\begin{equation} \label{eq:X-a-b-1}
\|X_{a,b}\|_{L_q} \leq \|a\|_{2} \Bigl(\|P_{J_1} b\|_2 + c \sqrt{\frac{q}{k}} \max_{\ell \geq 2} \|P_{J_\ell} b\|_2\Bigr).
\end{equation}

\subsection{Estimating $\Phi$}
Recall that
$$
\Phi=\sup_{t \in T} \sum_{s \geq s_1} \|A_\eps \Delta_s t\|_2
$$
and that for every $t \in \R^n$,
$$
Z_t = \sum_{i \not = j} \eps_i \eps_j t_i t_j \inr{Ae_i,Ae_j}.
$$

\begin{Theorem} \label{thm:Phi}
There are absolute constants $c$ and $C$ such that for $u >1$, with probability at least $1-2\exp(-cu^22^{s_1})$,
$$
\Phi^2 \leq C u^2 \delta^2 \ell_*^2(T).
$$
\end{Theorem}

\proof Let $s \geq s_1$, and as noted previously, for every $t \in T$,
\begin{align*}
\|A_\eps \Delta_s t\|_2^2 = & \sum_{i=1}^n \|Ae_i\|_2^2 (\Delta_s t)_i^2 + \sum_{i \not =j } \inr{Ae_i,Ae_j} \eps_i \eps_j (\Delta_st)_i (\Delta_st)_j
\\
\leq & (1+\delta)\|\Delta_s t\|_2^2 + \left| \sum_{i \not =j } \inr{Ae_i,Ae_j} \eps_i \eps_j (\Delta_st)_i (\Delta_st)_j\right|
\\
\leq & 2\|\Delta_s t\|_2^2+|Z_{\Delta_s t}|.
\end{align*}
Following the decoupling argument, one may fix $I \subset \{1,...,n\}$. The core of the argument is to obtain satisfactory estimates on moments of the random variables
\begin{equation*}
V_{\Delta_s t, I}= \sum_{j \in I^c} \eps_j (\Delta_s t)_j (W_{\Delta_st,I})_j
\end{equation*}
for that (arbitrary) fixed choice of $I$. And, since $|T_s| = 2^{2^s}$, for a uniform estimate that holds for every random variable of the form $V_{\Delta_s t,I}, \ t \in T$, it is enough to control the $q$-th moment of each $V_{\Delta_s t, I}$ for $q \sim 2^s$.

With that in mind, denote by $\E_{I^c}$ the expectation taken with respect to $(\eps_i)_{i \in I^c}$, and set $\E_{I}$ the expectation taken with respect to $(\eps_i)_{i \in I}$.

Let us apply \eqref{eq:X-a-b-1}, with the choice
$$
k=1/\delta^2 \leq k_*,  \ \ a=\Delta_s t \ \ {\rm and} \ \ b=W_{\Delta_st,I}.
$$
Thus, for every $(\eps_i)_{i \in I}$,
\begin{equation*}
\left(\E_{I^c} |V_{\Delta_s t, I}|^q \right)^{1/q} \leq \|\Delta_s t\|_{2} \left(\|P_{J_1} W_{\Delta_st,I}\|_2 + c_0 \delta \sqrt{q} \cdot \max_{\ell \geq 2} \|P_{J_\ell} W_{\Delta_st,I}\|_2 \right).
\end{equation*}
The sets $J_\ell$ are all of cardinality $1/\delta^2$ and so there are at most $\max\{\delta^2 n,1\}$ of them. By Lemma \ref{lemma:key-observation-1}, for each one of the sets $J_\ell \subset I^c$ and $q \geq |J_\ell| = 1/\delta^2$,
$$
\left(\E_I \|P_{J_\ell} W_{\Delta_st,I}\|_2^q \right)^{1/q} \leq c_1\delta \sqrt{q} \|\Delta_s t\|_2,
$$
implying that
\begin{align} \label{eq:max-est-1}
\left(\E_I \max_\ell \|P_{J_\ell} W_{\Delta_st,I}\|_2^q \right)^{1/q} \leq & \left(\sum_{\ell} \E_I \|P_{J_\ell} W_{\Delta_st,I}\|_2^q\right)^{1/q} \nonumber
\\
\leq & \ell^{1/q} \cdot c_1\delta \sqrt{q} \|\Delta_s t\|_2 \leq c_2\delta \sqrt{q} \|\Delta_s t\|_2
\end{align}
as long as $\max\{\delta^2 n,1\} \leq \exp(q)$. In particular, since $2^{s_1} \geq \log(e(1+n\delta^2))$, it suffices that $q \geq 2^{s_1}$ to ensure that \eqref{eq:max-est-1} holds.

Hence, for every $I \subset \{1,...,n\}$,
\begin{align*}
\left(\E |V_{\Delta_s t, I}|^q\right)^{1/q}  = & \left(\E_{I} \E_{I^c} |V_{\Delta_s t, I}|^q \right)^{1/q}
\\
\leq & \|\Delta_s t\|_{2} \left(\E_I \left(\|P_{J_1} W_{\Delta_st,I}\|_2 + c_0 \delta \sqrt{q} \cdot \max_{\ell \geq 2} \|P_{J_\ell} W_{\Delta_st,I}\|_2 \right)^q \right)^{1/q}
\\
\leq & \|\Delta_s t\|_{2} \cdot 2 \left( \left( \E_I \|P_{J_1} W_{\Delta_st,I}\|_2^q\right)^{1/q} + c_0 \delta \sqrt{q}\left(\E_I \max_{\ell \geq 2} \|P_{J_\ell} W_{\Delta_st,I}\|_2^q\right)^{1/q} \right)
\\
\leq & c_3\|\Delta_s t\|_{2}^2 \left(\delta\sqrt{q} + \delta^2q\right)
\\
\leq & c_4 \delta^2 q  \|\Delta_s t\|_{2}^2,
\end{align*}
because $\delta^2 q \geq \delta^2 2^{s_1} \geq 1$.

By Jensen's inequality, for every $t \in T$,
\begin{align*}
\left(\E|Z_{\Delta_s t}|^q\right)^{1/q} \leq & 4 \left(\E_\eta \E_\eps \left|\sum_{i \in I^c} \eps_i (\Delta_s t)_i (W_{\Delta_st,I})_i \right|^q\right)^{1/q}
\\
= & 4 \left(\E_\eta \E_\eps |V_{\Delta_s t,I}|^q \right)^{1/q} \leq c_5 q \delta^2 \|\Delta_s t\|_{2}^2.
\end{align*}
Setting $q=u2^s$ for $u \geq 1$ and $s \geq s_1$, Chebychev's inequality implies that with $(\eps_i)_{i=1}^n$-probability at least $1-2\exp(-c_6 u^2 2^s)$,
\begin{align*}
\|A_\eps \Delta_s t\|_2^2 = & 2\|\Delta_st\|_2^2+|Z_{\Delta_st}|  \leq  2\|\Delta_st\|_2^2+c_7 u^22^s \delta^2 \|\Delta_st\|_2^2
\\
\leq & c_8u^2\delta^22^s \|\Delta_st\|_2^2.
\end{align*}

By the union bound, this estimate holds for every $t \in T$ and $s \geq s_1$. Thus, there are absolute constants $c$ and $C$ such that with probability at least $1-2\exp(-c u^2 2^{s_1})$,
\begin{equation*}
\sum_{s \geq s_1} \|A_\eps \Delta_s t\|_2 \leq Cu\delta \sum_{s \geq s_1} 2^{s/2}\|\Delta_s t\|_2 \leq Cu\delta \ell_*(T),
\end{equation*}
as claimed.

\endproof

\subsection{Estimating $\Psi$}
Next, recall that
$$
\Psi^2=\sup_{t \in T} \left| \|A_\eps \pi_{s_1} t\|_2^2 - \|\pi_{s_1}t\|_2^2 \right|.
$$
Expanding as in \eqref{eq:expand}, the ``diagonal term" is at most $\delta^2 d_T^2$, and one has to deal with the ``off-diagonal" term
$$
\sup_{t \in T} |Z_{\pi_{s_1}t}|.
$$
As observed in Lemma \ref{lemma:key-observe-chaining} combined with Lemma \ref{rem:simple-Fubini}, it suffices to obtain, for every fixed $I \subset \{1,...,n\}$, estimates on the moments of
\begin{equation} \label{eq:Psi-est-1}
\sum_{j \in I^c} \eps_j (\Delta_s t)_j (W_{\pi_{s+1} t,I})_j \ \ \ \sum_{i \in I} \eps_i (\Delta_s t)_i (W_{\pi_{s} t,I^c})_i \ \ {\rm and} \ \   \sum_{j \in I^c} \eps_j (\pi_{s_0} t)_j (W_{\pi_{s_0} t,I})_j.
\end{equation}

\begin{Remark}
We shall assume throughout that $(\ell_*(T)/d_T)^2 \leq 1/\delta^2$; the required modifications when the reverse inequality holds are straightforward and are therefore omitted.
\end{Remark}

\vskip0.3cm

We begin with a standard observation:
\begin{Lemma} \label{lemma:moments}
There is an absolute constant $c$ such that the following holds. Let $(X_\ell)_{\ell \geq 0}$ be nonnegative random variables, let $q \geq 1$ and set $q_\ell=q2^\ell$. If there is some $\kappa$ such that for every $\ell$, $\|X_\ell\|_{L_{q_\ell}} \leq \kappa$, then
$$
\|\max_\ell X_\ell\|_{L_q} \leq c\kappa
$$
\end{Lemma}

\proof By Chebychev's inequality, for every $\ell \geq 1$ and $u \geq 2$,
$$
Pr\left( |X_\ell| \geq u\kappa\right) \leq \frac{\E|X_\ell|^{q_\ell}}{(u\kappa)^{q_\ell}} \leq u^{-q_\ell}.
$$
Therefore,
$$
Pr\left(\exists \ell \geq 2 : |X_\ell| \geq u\kappa \right) \leq \sum_{\ell \geq 2} u^{-q2^\ell} \leq c_1u^{-4q},
$$
implying that
$$
\E \max_{\ell \geq 2} |X_\ell|^q  \leq \int_0^\infty q u^{q-1} Pr\left(\max_{\ell \geq 1} |X_\ell| > u\right) du \leq (c_1 \kappa)^q.
$$
Hence,
$$
\|\max_\ell X_\ell\|_{L_q} \leq \|X_0\|_{L_q} + \|X_1\|_{L_q} + \|\max_{\ell \geq 2} X_\ell\|_{L_q} \leq c_2 \kappa.
$$
\endproof

\vskip0.3cm

The analysis is split into two cases.
\vskip0.3cm
\noindent {\bf Case I: $\log(e(1+\delta^2 n)) \leq 1/\delta^2$.}
\vskip0.3cm
In this case, $2^{s_1}=1/\delta^2$. For every $s_0 \leq s \leq s_1$ we invoke \eqref{eq:X-a-b-2} for sets $J_\ell$ of cardinality $2^{s+\ell}$ when $s+\ell \leq s_1$, and of cardinality $1/\delta^2$ when $s+\ell > s_1$.

Set $a=\Delta_s t$ and $b=W_{\pi_{s+1}t,I}$; the treatment when $b=W_{\pi_s t,I^c}$ is identical and is omitted. Let $q \geq 2^s$ and put $q_\ell=q2^\ell$ if $s+\ell \leq s_1$. Finally, set $p \geq 1/\delta^2$.

\vskip0.3cm
Consider $\ell$ such that $s+\ell \leq s_1$ and observe that $q_\ell \geq |J_\ell|=2^{s+\ell}$. By Lemma \ref{lemma:key-observation-1},
\begin{equation} \label{eq:Psi-eq-2}
|J_{\ell-1}|^{-1/2} \left(\E\|P_{J_\ell} W_{\pi_{s+1}t,I}\|_2^{q_\ell}\right)^{1/q_\ell} \leq c_1 \delta \|\pi_{s+1}t\|_2 \cdot \sqrt{\frac{q 2^\ell}{2^{s+\ell}}} \leq c_1 \sqrt{\frac{q}{2^s}}\delta d_T.
\end{equation}
Therefore, by Lemma \ref{lemma:moments},
$$
\left(\E \max_{\{\ell: s+\ell \leq s_1\}} \left(|J_{\ell-1}|^{-1/2} \|P_{J_\ell}W_{\pi_{s+1}t,I}\|_2\right)^q \right)^{1/q} \leq c_2 \sqrt{\frac{q}{2^s}}\delta d_T.
$$
Next, if $s+\ell > s_1$ then $p \geq |J_\ell|=1/\delta^2$ and by Lemma \ref{lemma:key-observation-1},
\begin{equation} \label{eq:Psi-eq-3}
|J_{\ell-1}|^{-1/2} \left(\E\|P_{J_\ell} W_{\pi_{s+1}t,I}\|_2^{p}\right)^{1/p} \leq c_3 \delta \|\pi_{s+1}t\|_2 \cdot \sqrt{\frac{p}{1/\delta^2}} \leq c_3 \sqrt{p\delta^2} \cdot \delta d_T.
\end{equation}
There are at most $n\delta^2$ sets $J_\ell$ in that range, and because $n\delta^2 \leq \exp(p)$, it is evident that
$$
\left(\E \max_{\{\ell: s+\ell > s_1\}} \left(|J_{\ell-1}|^{-1/2} \|P_{J_\ell} W_{\pi_{s+1}t,I}\|_2\right)^{p} \right)^{1/{p}} \leq c_4 \sqrt{p\delta^2}\cdot \delta d_T;
$$
therefore, as $q \leq p$,
$$
\left(\E \max_{\{\ell: s+\ell > s_1\}} \left(|J_{\ell-1}|^{-1/2} \|P_{J_\ell} W_{\pi_{s+1}t,I}\|_2\right)^{q} \right)^{1/{q}} \leq c_4 \sqrt{p\delta^2}\delta d_T;
$$
Thus, for every fixed $I \subset \{1,...,n\}$,
\begin{equation} \label{eq:max-over-sets}
\left(\E \max_{\ell} \left(|J_{\ell-1}|^{-1/2} \|P_{J_\ell} W_{\pi_{s+1}t,I}\|_2\right)^{q} \right)^{1/{q}} \leq c_5 \max\left\{\sqrt{\frac{q}{2^s}},\sqrt{p\delta^2}\right\} \cdot \delta d_T.
\end{equation}

Now, by \eqref{eq:X-a-b-1},
\begin{align} \label{eq:case1-in-proof}
& \left(\E_{I^c}\left| \sum_{j \in I^c} \eps_j (\Delta_s t)_j (W_{\pi_{s+1} t,I})_j\right|^q \right)^{1/q} \nonumber
\\
\leq &
\|\Delta_s t\|_{2} \left(\|P_{J_1} W_{\pi_{s+1}t,I}\|_2 + c_6 \sqrt{q} \max_{\ell \geq 2} \frac{\|P_{J_{\ell+1}} W_{\pi_{s+1}t,I}\|_2}{\sqrt{|J_{\ell}|}} \right),
\end{align}
which, combined with \eqref{eq:Psi-eq-2} and \eqref{eq:max-over-sets}, implies that
\begin{equation} \label{eq:case1-in-proof-a}
\left(\E \left| \sum_{j \in I^c} \eps_j (\Delta_s t)_j (W_{\pi_{s+1} t,I})_j\right|^q \right)^{1/q} \leq c_7 \sqrt{q} \|\Delta_s t\|_{2} \cdot \Bigl(1+ \max\left\{\sqrt{\frac{q}{2^s}},\sqrt{p\delta^2}\right\} \Bigr) \cdot \delta d_T.
\end{equation}
Clearly, there are at most $2^{2^{s+3}}$ random variables as in \eqref{eq:case1-in-proof-a}. With that in mind, set $u \geq 1$ and let $q=u2^s, \ p=u/\delta^2$. By Lemma \ref{lemma:key-observe-chaining} and Lemma \ref{rem:simple-Fubini} followed by the union bound, we have that with probability at least $1-2\exp(-c_8u^2 2^s)$, for every $t \in T$,
\begin{equation} \label{eq:case1-in-proof-b}
\E_\eta \left|\sum_{j \in I^c} \eps_j (\Delta_s t)_j W_{\pi_{s+1} t,I}\right| \leq c_9 u^2 2^{s/2} \|\Delta_s t\|_{2} \cdot \delta d_T.
\end{equation}
And, by the union bound and recalling that $2^{s_0}=(\ell_*(T)/d_T)^2$, \eqref{eq:case1-in-proof-b} holds uniformly for every $t \in T$ and $s_0 \leq s \leq s_1$ with probability at least
$$
1-2\exp\left(-c_{10}u^2 \left(\frac{\ell_*(T)}{d_T}\right)^2\right).
$$

An identical argument shows that with probability at least
$$
1-2\exp(-c_{11}u^2 2^{s_0})=1-2\exp(-c_{11}u^2 (\ell_*(T)/d_T)^2),
$$
for every $t \in T$,
$$
\E_\eta \left|\sum_{i \in I^c} \eps_i (\pi_{s_0} t)_i (W_{\pi_{s_0}t,I})_i\right| \leq c_{12} u^2 2^{s_0/2}  d_T \cdot \delta d_T \leq c_{13}d_T \delta \ell_*(T).
$$

Finally, invoking Lemma \ref{lemma:key-observe-chaining}, we have that with probability at least $1-2\exp(-cu^2(\ell_*(T)/d_T)^2)$, for every $t \in T_{s_1}$,
$$
|Z_t| \leq Cu^2 d_T \delta \ell_*(T),
$$
as required.
\endproof

\vskip0.3cm
\noindent {\bf Case II: $\log(e(1+\delta^2 n)) > 1/\delta^2$.}
\vskip0.3cm

The necessary modifications are minor and we shall only sketch them. In this range, $2^{s_1}=\log(e(1+\delta^2 n))$, and the problem is that for each vector $\Delta_s t$, the number of blocks $J_\ell$ of cardinality $1/\delta^2$---namely, $n\delta^2$, is larger than $\exp(2^s)$ when $s \leq s_1$. Therefore, setting $|J_\ell|=1/\delta^2$ for every $\ell$, the uniform estimate on
$$
\max_\ell \frac{\|P_{J_\ell+1} W_{\pi_{s+1}t,I}\|_2}{\sqrt{|J_{\ell}|}} =\delta \max_\ell \|P_{J_\ell} W_{\pi_{s+1}t,I}\|_2
$$
can be obtained by bounding $\left(\E \|P_{J_\ell} W_{\pi_{s+1}t,I}\|_2^q\right)^{1/q}$ for $q \geq \log (\delta^2 n)$. Indeed, by Lemma \ref{lemma:key-observation-1}, for every $I \subset \{1,...n\}$ and every $t \in T$, we have that
$$
\delta \left(\E \max_\ell \|P_{J_\ell} W_{\pi_{s+1}t,I}\|_2^q\right)^{1/q} \leq \delta^2 \sqrt{q} d_T.
$$
The rest of the argument is identical to the previous one and is omitted.
\endproof

\vskip0.3cm

\noindent{\bf Proof of Theorem \ref{thm:main-intro}.} Using the estimates we established, it follows that for $u \geq c_0$, with probability at least
$$
1-2\exp\left(-c_1u^2\left(\frac{\ell_*(T)}{d_T}\right)^2 \right),
$$
$$
\Phi \leq Cu \delta \ell_*(T) \ \ {\rm and} \ \ \Psi^2 \leq C u^2 d_T \delta \ell_*(T) \left(1+\delta \sqrt{\log(e(1+n\delta^2))}\right),
$$
noting that $\delta d_T^2 \leq  c_2d_T \delta \ell_*(T)$. Since
$$
\sup_{t \in T} \left| \|A_\eps t\|_2^2 - \|t\|_2^2 \right| \leq \Psi^2+2\Phi\sqrt{\Psi^2+d_T^2}+\Phi^2+C^\prime\left( (\delta \ell_*(T))^2 + d_T \delta \ell_*(T)\right)
$$
the claim follows from a straightforward computation, by separating to the cases $\Phi \leq \Psi$ and $\Psi \leq \Phi$.
\endproof

\section{Application - A circulant Bernoulli matrix} \label{sec:appl}
Let $(\eps_i)_{i=1}^n$ and $(\eps_i^\prime)_{i=1}^n$ be independent Bernoulli vectors. Set $\xi=(\eps_i^\prime)_{i=1}^n$ and put $D_\eps = {\rm diag}(\eps_1,...,\eps_n)$. Let $\Gamma$ be the circulant matrix generated by the random vector $\xi$; that is, $\Gamma$ is the matrix whose $j$-th row is the shifted vector $\tau_j\xi$, where for every $v \in \R^n$, $\tau_jv=(v_{j-i})_{i=1}^n$.

Fix $I \subset \{1,...,n\}$ of cardinality $m$ and let
$$
A=\sqrt{\frac{1}{m}}P_I \Gamma : \R^n \to \R^m
$$
to be the normalized partial circulant matrix. It follows from Theorem~3.1 in \cite{MR3263672} and the estimates in Section~4 there that for a typical realization of $\xi$, the matrix $A$ is $(\delta,k)$-regular for $\delta \sim m^{-1/2}\log^2n$:
\begin{Theorem} \cite{MR3263672} \label{thm:KMR}
There exist absolute constants $c$ and $C$ such that the following holds. For $x >0$ with probability at least
$$
1-2\exp\left(-c\min\left\{x^2\frac{m}{k},x\sqrt{\frac{m}{k}}\log^2 n \right\}\right)
$$
we have that
$$
\sup_{t \in \Sigma_k} \left| \|At\|_2^2 - 1 \right| \leq C(1+x)\sqrt{\frac{k}{m}}\log^2 n.
$$
\end{Theorem}
By Theorem \ref{thm:KMR} and the union bound for $1 \leq k \leq 1/\delta^2$, there is an event $\Omega$ with probability at least $1-2\exp(-c^\prime \log^4 n)$ with respect to $(\eps_i^\prime)_{i=1}^n$, on which, for every $k \leq 1/\delta^2$,
$$
\sup_{t \in \Sigma_k} \left| \|At\|_2^2 - 1 \right| \leq  \delta \sqrt{k}.
$$
This verifies the assumption needed in Theorem \ref{thm:main-intro} on the event $\Omega$. Now fix $(\eps_i^\prime)_{i=1}^n \in \Omega$ and let $A$ be the resulting partial circulant matrix. Set $A_\eps = A D_\eps$ and let $T \subset \R^n$. By Theorem \ref{thm:main-intro}, with probability at least
$$
1-2 \exp(-c^\prime (\ell_*(T)/d_T)^2)
$$
with respect to $(\eps_i)_{i=1}^n$, we have that
$$
\sup_{t \in T} \left| \|A_\eps t\|_2^2 - \|t\|_2^2 \right| \leq C^\prime \left( \Lambda d_T \frac{\ell_*(T)}{\sqrt{m}} \cdot \log^2 n + \left(\frac{\ell_*(T)}{\sqrt{m}}\right)^2 \cdot \log^4 n\right)
$$
where
$$
\Lambda \leq c^{\prime \prime} \max\left\{1,\frac{\log^5n}{m}\right\}.
$$

Thus, a random matrix generated by $2n$ independent random signs is a good embedding of an arbitrary subset of $\R^n$ with the same accuracy (up to logarithmic factors) as a gaussian matrix. Moreover, conditioned on the choice of the circulant matrix $A$, the probability estimate coincides with the estimate in the gaussian case.

\bibliographystyle{plain}
\bibliography{random}

\end{document}